\documentclass{amsart}
\usepackage{amsfonts,amsmath,amsthm,amssymb}
\usepackage[linesnumbered,ruled,noend]{algorithm2e}
\usepackage{comment}
\usepackage{float}
\usepackage{latexsym}
\usepackage{graphicx}
\usepackage{multicol}
\usepackage{color}
\usepackage{pifont}
\usepackage{enumerate}
\usepackage{lmodern}
\usepackage{makecell}
\usepackage{adjustbox}
\usepackage[backref,pdftex=true,colorlinks,linkcolor=blue,citecolor=red,pdfstartview=FitH]{hyperref}
\usepackage{setspace}
\usepackage[all]{xy}
\usepackage{pb-diagram,pb-xy}
\usepackage{pdflscape}
\usepackage{chngpage}
\usepackage{breakurl}
\usepackage[dvipsnames]{xcolor}

\makeatletter
\def\BState{\State\hskip-\ALG@thistlm}
\makeatother
\theoremstyle{plain}
\newtheorem{theorem}{Theorem}[section]

\theoremstyle{definition} 

\newtheorem{example}{Example}
\newtheorem*{example*}{Example}
\newtheorem{remark}[theorem]{Remark}

\newcommand{\Q}{\mathbb{Q}}
\newcommand{\Z}{\mathbb{Z}}

\begin{document}

\title[]{Computing rational solutions to $Ax^p+By^p+Cz^p=0$}

\author{Alejandro Arg\'{a}ez-Garc\'{i}a}
\address{Facultad de Matemáticas, Universidad Aut\'{o}noma de Yucat\'{a}n. Perif\'{e}rico Norte Kil\'{o}metro 33.5, Tablaje Catastral 13615 Chuburna de Hidalgo Inn, M\'{e}rida, Yucat\'{a}n, M\'{e}xico. C.P. 97200 }
\email{alejandro.argaez@correo.uady.mx}


\begin{abstract}
In this survey, we studied the possibility of finding rational solutions to the equation $Ax^p+By^p+Cz^p=0$ via its attached hyperelliptic curve $Y^2=X^p+A^2(BC)^{p-1}/4$ and its rational points computed using computational tools.
\end{abstract}
\maketitle

\section{Introduction}
The study of the Diophantine equation $Ax^p + By^p + Cz^p$ and its solutions
has been of interest for a long time, particularly in any number field. A classical approach is made by assuming there is a triplet $(x,y,z)$ such that $Ax^p + By^p + Cz^p=0$, then constructing either an elliptic curve or a hyperelliptic curve, depending on $p$, and classifying the points on that curve; for example, when $p=3$ and $K$ a quadratic field see \cite{ADP} and for $p\geq 5$ and $K$ a quadratic field see \cite{AP}.

In 1983, Faltings proved in \cite{faltings} that all curves of genus at least two, over any number field $K$, have finitely many $K$-rational points. Thus, we know there are only a finite number of triplets where the rational points of the hyperelliptic curve can come. While Falting's theorem tells us that there are finitely many points, it does not tells us how many are. In order to try to count how many points there are, we use the Chabauty-Coleman's bound \cite{chabauty}, \cite{coleman}, which provide us with an idea of how many points could be when the rank of the Jacobian is less than the genus of the curve.

In this survey, we explore how to compute rational solutions to $Ax^p+By^p+Cz^p=0$ by calculating all rational points of its associated hyperelliptic curve over $\Q$ using computational tools, specifically, using Magma's online calculator.
 
\section{Diophantine equations and hyperelliptic curves} 
Let $(x,y,z)$ be any solution to the Diophantine equation 
\begin{equation}\label{eq:00}
Ax^p+By^p+Cz^p=0
\end{equation}
for $p\geq 5$ prime, $x\neq 0$, and $A$, $B$ and $C$ being pairwise coprime $p$th powerfree integers in $\Z$. There is a standard change of variable (Proposition 6.4.13. \cite{cohen01}) to obtain the hyperelliptic curve 

\begin{equation}\label{eq:01}
E\colon Y^2 = X^{p}+\dfrac{A^2(BC)^{p-1}}{4}
\end{equation}
given by 
\begin{equation}\label{eq:cambiodevariable} 
X=\dfrac{-BCyz}{x^2}, \quad Y=\dfrac{(-BC)^{\frac{(p-1)}{2}}(By^p-Cz^p)}{(2x^{p})}
\end{equation}

In other words, all triplets $(x,y,z)$ that are a solution to $(\ref{eq:00})$ are sent to points $(X,Y)$ in $E(\Q)$ given by $(\ref{eq:cambiodevariable})$, but not all points in $E(\Q)$ necessarily have an associated triplet $(x,y,z)$ satisfying $(\ref{eq:00})$. Thus, in order to determine whether there are triplets, we need to be able to compute all rational points in $E(\Q)$ and construct a method to find the triplet associated with the point we chose. 

Finally, recall that the degree of the polynomial defining the hyperelliptic curves is the prime $p=2g+1$, where $g$ is the genus of the curve $E$, so we only have one point at infinity, named $(1,0,0)$.

\section{Ordering $A$, $B$ and $C$}

How we order $A$, $B$, and $C$ matters, specifically in the construction of the hyperelliptic curve and its rational points. Taking the ordering $(A,*,*)$, we have eight possible generalized Fermat equations 
\[
\pm Ax^p\pm By^p\pm Cz^p=0\quad \text{and}\quad \pm Ax^p\pm Cy^p\pm Bz^p=0
\]
where all their solutions go to the same hyperelliptic curve 
\[
Y^2=X^p+\dfrac{A^2(BC)^{p-1}}{4}
\]
where all their solutions go to the same hyperelliptic curve. On the other hand, taking the ordering $(B, *, *)$ give us another eight possible generalized Fermat equations 
\[
\pm Bx^p\pm Ay^p\pm Cz^p=0\quad\text{and}\quad \pm Bx^p\pm Cy^p\pm Az^p=0
\]
where all their solutions go to the same hyperelliptic curve 
\[
Y^2=X^p+\dfrac{B^2(AC)^{p-1}}{4}
\]
and the last ordering $(C,*,*)$ gives also eight possible generalized Fermat equations 
\[
\pm Cx^p\pm Ay^p\pm Bz^p=0\quad\text{and}\quad \pm Cx^p\pm By^p\pm Az^p=0
\]
where all their solutions go to the same hyperelliptic curve
\[
Y^2 = X^p+\dfrac{C^2(AB)^{p-1}}{4}
\]
We need to be vigilant in the ordering of each triplet, but we do not need to worry about choosing their signs; the following theorem, which is a simple observation, confirms it.

\begin{theorem}\label{thm:signdoesntmatter}
Let $(x,y,z)$ be a solution to $Ax^p+By^p+Cz^p=0$, then the following are equivalent
\begin{enumerate}[$(a)$]
\item $(-x,y,z)$ is a solution to $-Ax^p+By^p+Cz^p=0$
\item $(x,-y,z)$ is a solution to $Ax^p-By^p+Cz^p=0$
\item $(x,y,-z)$ is a solution to $Ax^p+By^p-Cz^p=0$
\item $(-x,-y,z)$ is a solution to $-Ax^p-By^p+Cz^p=0$
\item $(-x,y,-z)$ is a solution to $-Ax^p+By^p-Cz^p=0$
\item $(x,-y,-z)$ is a solution to $Ax^p-By^p-Cz^p=0$
\item $(-x,-y,-z)$ is a solution to $-Ax^p-By^p-Cz^p=0$
\end{enumerate}
\end{theorem}

Therefore, we can focus only on finding rational solutions to the equations 
\begin{align*}
Ax^p+By^p+Cz^p&=0\\
Bx^p+Ay^p+Cz^p&=0\\
Cx^p+By^p+Az^p&=0
\end{align*}
for the set $\{A,B,C\}$. We end this section with the following theorem, which again is another simple observation and a handy one. 

\begin{theorem}\label{thm:fact}
If $A+B+C=0$, then $(1,1,1)$ is a solution to $Ax^p+By^p+Cz^p=0$.
\end{theorem}

\section{Finding solutions to $Ax^p+By^p+Cz^p=0$}

We want to compute the triplet $(x,y,z)$ such that $Ax^p+By^p+Cz^p=0$. To do this, take a rational point $(X,Y)$ on $(\ref{eq:01})$ satisfying $XY\neq 0$ and take
\[
\pm Y=\dfrac{(-BC)^{\frac{(p-1)}{2}}(By^p-Cz^p)}{(2x^{p})}
\]
meaning $(x,y,z)$ satisfies the equation
\[
\pm\dfrac{2Y}{(-BC)^{\frac{(p-1)}{2}}}x^p-By^p+Cz^p=0
\]
In other words, $(x,y,z)$ is the solution of the system of equations
\begin{align*}
Ax^p+By^p+Cz^p&=0\\
\pm A'x^p-By^p+Cz^p&=0
\end{align*}
where $A'=\dfrac{2Y}{(-BC)^{\frac{(p-1)}{2}}}$, which a priori is not necessarily an integer number. By linear algebra, we have 
\[
\begin{pmatrix}
A & B\\
\pm A' & -B
\end{pmatrix}
\begin{pmatrix}
x^p\\
y^p
\end{pmatrix}
=
\begin{pmatrix}
-Cz^p\\
-Cz^p
\end{pmatrix}
\]
Thus, we have
\begin{equation}\label{eq:beforesqrtsolution}
\dfrac{x^p}{z^p}=-\dfrac{2C}{(A\pm A')} \quad \dfrac{y^p}{z^p}=\dfrac{(-A\pm A')C}{(-A\mp A')B}
\end{equation}
and therefore 
\begin{equation}\label{eq:sqrtpsolution}
\sqrt[p]{-\dfrac{2C}{(A\pm A')}}\quad \text{and} \quad \sqrt[p]{\dfrac{(-A\pm A')C}{(-A\mp A')B}}
\end{equation}
If these number were rational numbers, then we would have found the triplet $(x,y,z)$.

\begin{remark}
We have to choose $XY\neq 0$, because if not, by $(\ref{eq:cambiodevariable})$,  for $X=0$ we would have either $y=0$ or $z=0$, and for $Y=0$ we would have $\sqrt[p]{B/C}\in\Q$, which is a contradiction.
\end{remark}

\section{Examples}
All the computations were done in Magma's online calculator \cite{magmacalculator} and the setting \texttt{SetClassGroupBounds("GRH");}.

\begin{example}
Let $\{121,123,125\}$ be a set of coprime integer numbers and consider the generalized Fermat equation $$123x^5+125y^5+121z^5=0$$ Its attached hyperelliptic curve is 
\begin{equation}\label{example:01}
Y^2 =X^5+\dfrac{121^2 123^4 125^4}{4}\Leftrightarrow N^2 = M^5+202689719415562500000000
\end{equation}
We look for rational points on $(\ref{example:01})$ using \texttt{Points} in Magma, but we are only able to obtain 
\[
\{(1:0:0),(0:-450210750000:1),(0:450210750000:1) \}
\]
Nevertheless, its Jacobian has rank $0$, which was computed using \texttt{JacobianRank}, meaning we can use the implementation \texttt{Chabauty0} and confirm that those are the only rational points to our hyperelliptic curve. Therefore, there are no rational triplets $(x,y,z)$ satisfying $123x^5+125y^5+121z^5=0$.
\end{example}

\begin{example}
Let $\{2,9,11\}$ be a set of coprime integer numbers and consider the generalized Fermat equation $$2x^5+9y^5+11z^5=0$$ It is straightforward that the triplets $(1,1,-1)$ and $(-1,-1,1)$ are solutions. We want to find more triplets $(x,y,z)$ such that $2x^5+9y^5+11z^5=0$ with $xyz\neq \pm 1$.

First, we construct the hyperelliptic curve
\begin{equation}\label{example:2x9y11z}
Y^2=X^5+96059601
\end{equation}

Its Jacobian rank is $1$, meaning we can use the \texttt{Chabauty} implementation on Magma to get all rational points on $(\ref{example:2x9y11z})$, which are
\[
\{(1:0:0),(0:-9801:1),(0:9801:1),(99:-98010:1),(99:98010:1)\}
\]
By applying $(\ref{eq:cambiodevariable})$ and $(\ref{eq:beforesqrtsolution})$ to those points with $XY\neq 0$, we get $x^2=-yz$ and
\[
\dfrac{x^5}{z^5} = -1, \dfrac{y^5}{z^5} = -1
\]
which means $x=-z$ and $y=-z$. In this way, the triplet $(-z,-z,z)$ give us all rational solutions. On the other hand, we have $x^2=-yz$ and 
\[
\dfrac{x^5}{z^5} = \dfrac{11}{9}, \dfrac{y^5}{z^5} = \dfrac{-121}{81}
\]
meaning
\[
\sqrt[5]{\dfrac{11}{9}},\sqrt[5]{\dfrac{-121}{81}}\not\in\Q
\]
Therefore, there are no other rational solutions.
\end{example}

\begin{example}
Let $\{16,9,7\}$ be a set of coprime integer numbers and consider the generalized Fermat equation $$16x^5+9y^5+7z^5=0$$ It is straightforward that the triplets $(1,-1,-1)$ and $(-1,1,1)$ are solutions. We want to find more triplets $(x,y,z)$ such that $16x^5+9y^5+7z^5=0$ with $xyz\neq \pm 1$.

First, we construct the hyperelliptic curve
\begin{equation}\label{example:16x9y7z}
Y^2=X^5+1008189504
\end{equation}

Its Jacobian rank is $1$, meaning we can use the \texttt{Chabauty} implementation on Magma to get all rational points on $(\ref{example:16x9y7z})$, which are
\[
\{(1:0:0),(0:-31752:1),(0:31752:1),(-63:-3969:1),(-63:3969:1)
\]
By applying $(\ref{eq:cambiodevariable})$ and $(\ref{eq:beforesqrtsolution})$ to those points with $XY\neq 0$, we get $x^2=-yz$ and
\[
\dfrac{x^5}{z^5} = -1, \dfrac{y^5}{z^5} = 1
\]
which means $x=-z$ and $y=z$. In this way, the triplet $(-z,z,z)$ give us all rational solutions. On the other hand, we have $x^2=-yz$ and 
\[
\dfrac{x^5}{z^5} = \dfrac{-7}{9}, \dfrac{y^5}{z^5} = \dfrac{49}{81}
\]
meaning
\[
\sqrt[5]{\dfrac{-7}{9}},\sqrt[5]{\dfrac{49}{81}}\not\in\Q
\]
Therefore, there are no other rational solutions.
\end{example}

\end{document}